# THE HANNA NEUMANN CONJECTURE IS TRUE WHEN ONE SUBGROUP HAS A POSITIVE GENERATING SET

Bilal Khan

**Abstract**

The **Hanna Neumann conjecture** states that if $F$ is a free group, then for all finitely generated subgroups $H, K \leqslant F$,

$$\mathrm{rank}(H \cap K) - 1 \leqslant [\mathrm{rank}(H) - 1]\,[\mathrm{rank}(K) - 1]$$

In this paper, we show that if one of the subgroups, say $H$, has a generating set consisting of only positive words, then $H$ is not part of any counterexample to the conjecture. We further show that if $H \leq F_2$ is part of a counterexample to the conjecture, then its folding $\Gamma_H$ must contain source and/or sink vertices.

## 1  Introduction

Improving Howson's earlier bound [4] on the rank of intersections of finitely generated (f.g.) subgroups of free groups, H. Neumann proved in [7] that any $H, K \leqslant_{\mathrm{f.g.}} F$ must satisfy

$$rank(H \cap K) - 1 \leq 2[rank(H) - 1][rank(K) - 1]$$

The stronger assertion obtained by omitting the factor of 2 has come to be known as the Hanna Neumann conjecture. In [1], Burns improved H. Neumann's bound by showing that in fact

$$rank(H \cap K) - 1 \leq 2[rank(H) - 1][rank(K) - 1] - min(rank(H) - 1, rank(K) - 1)$$

In 1983, J. Stallings introduced the notion of a *folding* and showed how to apply these objects in the study of subgroups of free groups [11]: Recall the well-defined constructive map which assigns to each finitely generated subgroup $H$ of a free group $F = F(X)$, a corresponding folding $\Gamma_H = (V_H, E_H)$. We view the folding $\Gamma_H$ as a deterministic finite automaton, represented as a directed multigraph (with loops), with each directed edge labeled by an element of the ground set $X$. $\Gamma_H$ enjoys the property that the set of freely reduced elements in $H$ coincides with the set of words that can be read along closed non-backtracking walks that start and end at a distinguished vertex $1_H \in V_H$.

Stallings's approach was applied by Gersten in [3] to solve certain special cases of the conjecture. Similar techniques were developed over a sequence of papers by Imrich [6, 5],

Nickolas [9], and Servatius [10] who gave alternate proofs of Burns' bound and resolved special cases of the conjecture. In 1989, W. Neumann showed that the conjecture is true "with probability 1" for randomly chosen subgroups of free groups [8], and proposed a stronger form of the conjecture. In 1992, Tardos proved in [12] that the conjecture is true if one of the two subgroups has rank 2. In 1994, Warren Dicks showed that the strong Hanna Neumann conjecture is equivalent to a conjecture on bipartite graphs, which he termed the Amalgamated Graph conjecture [2]. In 1996, Tardos used Dicks' method to give the first new bound for the general case in [13], where he proved that $\forall H, K \leqslant_{\text{f.g.}} F$,

$$rank(H \cap K) - 1 \leq 2[rank(H) - 1][rank(K) - 1] - rank(H) - rank(K) + 1$$

To date, this is the best bound for the general case; the conjecture remains open.

Research in this area has focused largely on "translating" the group-theoretic properties of subgroups of free groups into the graph-theoretic properties of their corresponding foldings. In 1999, at the NY Group Theory Seminar, A. Miasnikov proposed a research project to elaborate the reverse, i.e. to interpret well-known properties of graphs in group-theoretic terms. In this paper we present the results of this ongoing project, and present some of the consequences for the conjecture.

The outline of the paper is as follows. We introduce the notion of a strong directed trail decomposition for directed graphs, and demonstrate that the existence of such a decomposition is equivalent to strong connectivity of the directed graph. Then, we show that strong connectivity of the folding of a subgroup $H \leqslant_{\text{f.g.}} F$ is equivalent to the existence of a set of positive words that generate $H$. Finally, we show that if a folding of $H \leqslant_{\text{f.g.}} F_2$ has a directed trail decomposition, then the folding necessarily exhibits a symmetry in the distribution of its degree 3 vertices; we term this the "3-balanced" property. Invoking a result of W. Neumann [8] we argue that a group $H \leqslant_{\text{f.g.}} F_2$ which has a 3-balanced folding cannot be part of any counterexample to the Hanna Neumann conjecture. Thus, we prove that every pair $H, K \leqslant_{\text{f.g.}} F_2$ satisfies the Hanna Neumann conjecture when least one of the two groups is generated by a finite set of positive words. By an embedding argument, we show the same result holds for finitely generated subgroups of any free group $F$.

## 2 Definitions

For concreteness, much of this exposition is restricted to the finitely generated subgroups of $F_2 = F(\{a, b\})$. Suppose we are given $H$, a non-trivial finitely generated subgroup of $F_2$, and its folding $\Gamma_H = (V_H, E_H)$. The reader who wishes to review a standard constructive definition of $\Gamma_H$ may consult the proof of lemma 3.4 on page 7, where it is outlined. Since $H \leqslant_{\text{f.g.}} F_2$, $\Gamma_H$ has vertices of undirected degree $\leqslant 4$ (the undirected degree of a vertex is the sum of its in-degree and out-degree). Figure 1 illustrates the types of vertices that may be present in $\Gamma_H$. Define $d = d_H : V_H \to \{1, 2, 3, 4\}$ to be the function that assigns to each vertex $v \in V_H$ its undirected degree in $\Gamma_H$. Now put $d_i(\Gamma_H) = |\{v \in V_H|\ d_H(v) = i\}|$, for $i = 1, 2, 3, 4$. We classify vertices of degree 3 based on the labels of their incident edges, naming the 4 classes $C_1, C_2, C_3$, and $C_4$; these classes are shown in figure 1. Define $C_i(\Gamma_H)$ to be the number of degree 3 vertices of type $C_i$ in $\Gamma_H$.

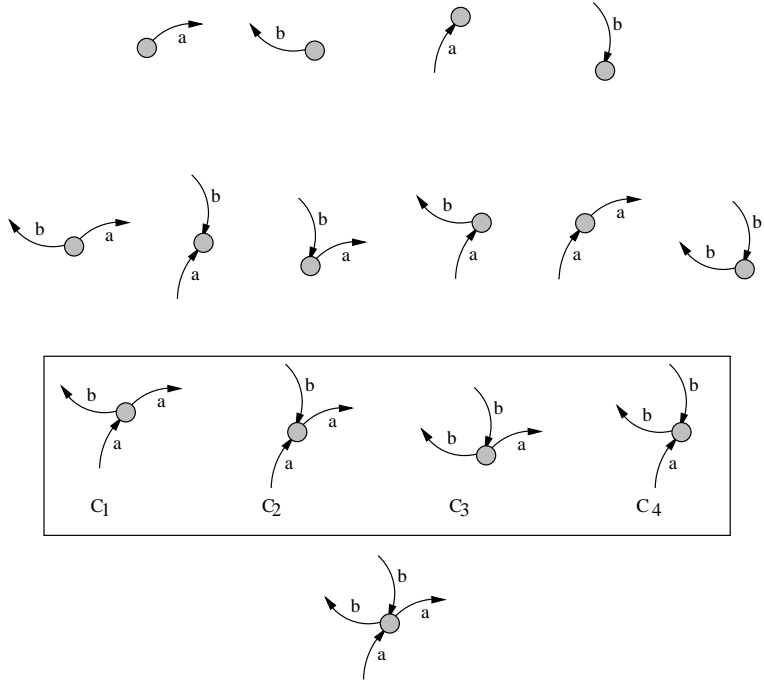

Figure 1: Local structure of vertex labels in a folding.

**Definition 2.1.** A folding $\Gamma$ is called **3-balanced** if it satisfies

$$C_1(\Gamma) + C_3(\Gamma) = C_2(\Gamma) + C_4(\Gamma) \tag{1}$$

**Remark 2.2.** In [8], Walter Neumann showed that if $H, K \leqslant_{\text{f.g.}} F_2$ are a counterexample to the conjecture, then $\exists i \in \{1, 2, 3, 4\}$ s.t. $C_i(\Gamma_H) > \frac{1}{2}d_3(\Gamma_H)$ and $C_i(\Gamma_K) > \frac{1}{2}d_3(\Gamma_K)$. Clearly, if a group has a 3-balanced folding, then no more than half of its degree 3 vertices can be of the same type. Thus, it follows from W. Neumann's result that if $H$ has a 3-balanced folding, then there is no $K \leqslant_{\text{f.g.}} F_2$ for which the pair $(H, K)$ are a counterexample to the conjecture.

**Definition 2.3.** A **directed trail** $P$ in a directed graph $\Gamma = (V, E)$ is a non-empty sequence of distinct directed edges $e_1, e_2, \ldots, e_m$ ($e_i \in E$, $i = 1, \ldots, m$) for which $tail(e_{i+1}) = head(e_i)$ ($i = 1, \ldots, m-1$). The length of $P$ is denoted $|P| = m$. The *start* of $P$ is denoted as $s(P) = tail(e_1)$, and the *terminus* of $P$ as $t(P) = head(e_{|P|})$.

**Definition 2.4.** A **self-avoiding directed trail** $P = (e_1, e_2, \ldots, e_m)$ in a directed graph $\Gamma = (V, E)$ is a directed trail which additionally satisfies $\forall i, j \in \{1, \ldots, m\}$,

$$i \neq j \Rightarrow [tail(e_i) \neq tail(e_j) \text{ and } head(e_i) \neq head(e_j)]$$

Notice that by this definition, a self-avoiding trail may satisfy $head(e_m) = tail(e_1)$.

**Remark 2.5.** Given a directed trail $P = (e_1, e_2, \ldots, e_{|P|})$ in a directed graph $\Gamma = (V, E)$, it is easy to verify that there always exists a *self-avoiding* trail $P' = (f_1, f_2, \ldots, f_{|P'|})$, where for $i = 1, \ldots, |P|$, $f_i \in \{e_1, e_2, \ldots, e_{|P|}\}$ and $s(P) = s(P')$, $t(P) = t(P')$.

The following construction is inspired by the "open ear decomposition" of Whitney [14].

**Definition 2.6.** A **directed trail decomposition** of a directed graph $\Gamma = (V, E)$ is a sequence of directed trails $\mathcal{P} = (P_0, \ldots, P_n)$ satisfying the following 3 conditions:

1. The trails are a partition of the edges of $\Gamma$:
$$\bigcup_{i=0}^{n} P_i = E \text{ and } i \neq j \Rightarrow P_i \cap P_j = \emptyset$$

2. $s(P_0) = t(P_0)$, and we denote this vertex as $1_\Gamma$.

3. For each $i = 1, \ldots, n$, the directed trail $P_i$ satisfies:
$$V[P_i] \cap \bigcup_{j=0}^{i-1} V[P_j] \neq \emptyset \quad \Rightarrow \quad V[P_i] \cap \bigcup_{j=0}^{i-1} V[P_j] = \{s(P_i), t(P_i)\}$$
$$V[P_i] \cap \bigcup_{j=0}^{i-1} = \emptyset \quad \Rightarrow \quad s(P_i) = t(P_i)$$

**Definition 2.7.** A **strong directed trail decomposition** of a directed graph $\Gamma = (V, E)$ is a directed trail decomposition $\mathcal{P}$ of $\Gamma$ which satisfies
$$\forall i \in \{1, \ldots, |\mathcal{P}|\}, \ V[P_i] \cap \bigcup_{j=0}^{i-1} V[P_j] \neq \emptyset$$

**Definition 2.8.** Given a directed graph $\Gamma = (V, E)$ define the binary relation $SC \subseteq V \times V$ (**strong connectivity**). Specifically, for $u, v \in V$

$$(u, v) \in SC \leftrightarrow \begin{array}{l} u = v, \text{ or} \\ \text{[There is a directed trail from } u \text{ to } v, \\ \text{and there is a directed trail from } v \text{ to } u] \end{array}$$

Let $J_1, \ldots J_m$ denote the equivalence classes of $V$ under $SC$, and define the $i$th **strongly connected component** $\widetilde{J_i}$ to be the subgraph induced by $J_i$ ($i = 1, \ldots, m$). A directed graph $\Gamma = (V, E)$ is **strongly connected** iff $SC = V \times V$.

**Definition 2.9.** A strongly connected component $\widetilde{J}$ is referred to as a **source** [resp. **sink**] if $\widetilde{J}$ consists of a single vertex $v$, where $v$ has degree 2 in $\Gamma$, and $v$ is incident to exactly two outgoing [resp. incoming] edges. A group $H \leqslant_{\text{f.g.}} F(X)$ is called **source/sink-free** (with respect to basis X) if $\Gamma_H$ contains neither source nor sink vertices.

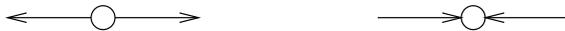

Figure 2: A source vertex (left) and a sink vertex (right).

The property of being source/sink-free depends on the choice of basis for the ambient free group $F$. This leads us to define the following basis-invariant version of this property:

`Definition 2.10.` A subgroup $H \leqslant_{\text{f.g.}} F$ is termed **potentially source/sink-free** if $\exists \phi \in Aut(F)$, such that $\Gamma_{\phi(H)}$ is source/sink-free (with respect to basis X).

`Definition 2.11.` Given a word $w \in F(X) = F(\{x_1, \ldots, x_n\})$, we say that $w$ is **positive** (with respect to basis $X$) if its freely reduced form consists only of the symbols $x_1, \ldots, x_n$ (and contains no occurrence of $x_1^{-1}, \ldots, x_n^{-1}$). A word $w$ is called **negative** if and only if $w^{-1}$ is a positive word (with respect to basis $X$). A subgroup $H \leqslant F(X)$ is said to have a positive generating set (or simply: "H is positively generated") if $\exists S \subseteq H$ such that $\langle S \rangle = H$ and $\forall w \in S$, $w$ is positive. Note that unless explicitly stated, $S$ need not be a basis.

The property of being positively generated depends on the choice of basis for the ambient free group $F$. This leads us to define the following basis-invariant version of this property:

`Definition 2.12.` A subgroup $H \leqslant_{\text{f.g.}} F(X)$ is called **potentially positive** if for some $\phi \in Aut(F)$, $\phi(H)$ is positively generated (with respect to basis X).

## 3  Results

The results of this paper are organized as follows. First, in lemmas 3.1 and 3.2, we show that the strong connectivity of a folding is equivalent to the existence of a strong directed trail decomposition. Then, lemma 3.4 demonstrates that if a finitely generated subgroup of a free group has a generating set consisting of only positive words, then the group's associated folding is necessarily strongly connected. Lemma 3.5 proves the converse: that strong connectivity of a folding implies the existence of a positive basis for its associated group. Thus, along the way, we determine (in corollary 3.6) that the existence of a positive generating set for a subgroup $H \leqslant_{\text{f.g.}} F$ is equivalent to the existence of a positive basis for $H$. Finally, in lemma 3.7, we show that if a folding has a directed trail decomposition, then it is necessarily 3-balanced. These results are combined in theorem 3.9 to show that if $H, K$ are two finitely generated subgroups of a free group $F$ and at least one of the two subgroups is generated by a set of positive words, then the pair $(H, K)$ satisfy the Hanna Neumann conjecture.

`Lemma 3.1.` *If $\Gamma = (V, E)$ has a strong directed trail decomposition, then $\Gamma$ is a strongly connected directed graph.*

*Proof.* Let $P_0, \ldots, P_n$ be a strong directed trail decomposition of $\Gamma$. Let $v_0 \in V$ be arbitrary; we show there is a directed trail from $v_0$ to $1_\Gamma$ and a directed trail from $1_\Gamma$ to $v_0$.

Clearly, $v_0 \in P_{i_0}$ for some $i_0 \in \{1, \ldots, n\}$. Successively, for each $m \geqslant 1$, if $i_{m-1} \neq 0$ we define $v_m = t(P_{i_{m-1}})$ and choose $i_m < i_{m-1}$ such that $v_m \in P_{i_m}$. Since $i_0, i_1, i_2, \ldots, i_m, \ldots$ are monotonically decreasing indices from the finite set $\{0, \ldots n\}$, there is some $M$ sufficiently large for which $i_M = 0$, and hence $v_M = 1_\Gamma$. By concatenating final segments of the directed trails $P_{i_0}, P_{i_1}, \ldots, P_{i_{M-1}}$, we obtain the directed trail

$$v_0 \stackrel{P_{i_0}}{\rightsquigarrow} v_1 \stackrel{P_{i_1}}{\rightsquigarrow} \ldots \stackrel{P_{i_{M-1}}}{\rightsquigarrow} v_M$$

which connects $v_0$ to $1_\Gamma$.

Put $j_0 = i_0$. Successively, for each $\ell \geqslant 1$, if $j_{\ell-1} \neq 0$ we define $u_\ell = s(P_{j_{\ell-1}})$ and choose $j_\ell < j_{\ell-1}$ such that $u_\ell \in P_{j_\ell}$. Since $j_0, j_1, j_2, \ldots, j_\ell, \ldots$ are monotonically decreasing indices from the finite set $\{0, \ldots n\}$, there is some $L$ sufficiently large for which $j_L = 0$, and hence $u_L = 1_\Gamma$. By concatenating initial segments of the directed trails $P_{j_{L-1}}, P_{j_{L-2}}, \ldots, P_{j_0}$, we obtain the directed trail

$$u_L \stackrel{P_{j_{L-1}}}{\rightsquigarrow} u_{L-1} \stackrel{P_{j_{L-2}}}{\rightsquigarrow} \ldots \stackrel{P_{j_1}}{\rightsquigarrow} u_1 \stackrel{P_{j_0}}{\rightsquigarrow} v_0$$

connecting $1_\Gamma$ to $v_0$. □

The converse of lemma 3.1 is also true, as we now show.

**Lemma 3.2.** *If $\Gamma = (V, E)$ is a strongly connected directed graph, then $\Gamma$ has a strong directed trail decomposition consisting of self-avoiding directed trails.*

*Proof.* Fix an arbitrary vertex in $\Gamma$, and denote it as $1_\Gamma$. We give the following effective procedure for constructing a directed trail decomposition. First, define $\Gamma_0 \stackrel{\text{def}}{=} (1_G, \emptyset)$. Then, starting with $i = 0$:

1. If $V[\Gamma_i] = V$, proceed to step 2. Otherwise, fix any $v \in V \backslash V[\Gamma_i]$. Since $\Gamma$ is strongly connected, we can choose a directed trail $s$ from $1_\Gamma$ to $v$, and a directed trail $t$ from $v$ to $1_\Gamma$. Suppose that $s$ is the sequence of edges $(s_1, s_2, \ldots, s_{|s|})$, where $tail(s_1) = 1_\Gamma$ and $head(s_{|s|}) = v$. Fix $s_j$ to be the last edge in $s$ for which $tail(s_j) \in V[\Gamma_i]$; put $x = tail(s_j)$. Similarly, suppose that $t$ is the sequence of edges $(t_1, \ldots, t_{|t|})$, where $tail(t_1) = v$ and $head(t_{|t|}) = 1_\Gamma$. Fix $t_k$ to be the first edge in $t$ for which $head(t_k) \in V[\Gamma_i]$; put $y = head(t_k)$. Define

$$P_i = (s_j, s_{j+1}, \ldots, s_{|s|}, t_1, \ldots t_{k-1}, t_k)$$

   Clearly $P_i$ is a trail. In light of remark 2.5, we may by suitably adjusting our choice of $v$, assume that $P_i$ is a self-avoiding trail. Put $\Gamma_{i+1} = \Gamma_i \sqcup P_i$. Increment $i$. Repeat step 1.

2. If $E[\Gamma_i] = E$, halt. Otherwise, fix any $e = (x, y) \in E \backslash E[\Gamma_i]$. Take the $i$th directed trail to be $P_i = (e)$ and put $\Gamma_{i+1} = \Gamma_i \sqcup P_i$. Increment $i$. Repeat step 2.

Notice that at each iteration of the procedure, the trail $P_i$ is constructed so that it does not contain any edges already in $\Gamma_i$. Indeed, for $i \geq 0$, $P_i$ attaches to $\Gamma_i$ at precisely its start and terminus vertices $x, y$. Thus, the procedure outputs a strong directed trail decomposition of $\Gamma$ that consists of self-avoiding directed trails. □

We shall later need the following technical refinement of the above lemma.

**Corollary 3.3.** *If $\Gamma = (V, E)$ is a strongly connected directed graph, then $\Gamma$ has a strong directed trail decomposition consisting of self-avoiding directed trails $P_0, \ldots, P_n$, which additionally satisfy the condition that for each $i = 1, \ldots n+1$, $\Gamma_i = \bigcup_{j=0}^{i-1} P_j$ is a strongly connected directed graph.*

*Proof.* Note that in the proof of lemma 3.2, for each $i = 1, \ldots n + 1$, the sequence $P_0, \ldots, P_{i-1}$ forms a strong directed trail decomposition of $\Gamma_i$. Thus, by lemma 3.1, $\Gamma_i$ is a strongly connected directed graph. □

The next lemma demonstrates that if a subgroup $H$ of a free group $F$ is generated by a positive set of words, then $H$ must necessarily have a strongly connected folding.

**Lemma 3.4.** *Let $H$ be a finitely generated subgroup of a free group $F(X)$, and let $\Gamma = (V, E)$ be the folding of $H$. If $H$ is positively generated (with respect to basis $X$), then $\Gamma_H$ is a strongly connected directed graph.*

*Proof.* Let $S = \{w_1, \cdots, w_n\}$ be a positive generating set for $H$, i.e. $H = \langle S \rangle$, where $w_i$ is positive for all $i = 1, \ldots, n$.

Let $R_H$ be the graph obtained as follows: (1) Construct $n$ directed cycles $c_1 = (V_1, E_1)$, ..., $c_n = (V_n, E_n)$, where $|V_i| = |w_i|$. (2) Pick one vertex from each of the cycles, and identify this subset of vertices; denote the resulting vertex $1_H$. (3) Label cycle $c_i$'s edges by successive letters of $w_i$. We call $R_H$ the *rose* of $H$. Because the generating set $S$ consists of positive words, $R_H$ contains a directed path $p_{u,v}$ between any two vertices $u$ and $v$—simply take $p_{u,v}$ to be the path that goes from $u$ to the vertex $1_H$ followed by the path from $1_H$ to $v$ (see figure 3).

Now recall the "folding process" $\pi$ that transforms $R_H$ into the folding $\Gamma_H$: Repeatedly identify pairs of edges $e, e'$ which satisfy

$$label(e) = label(e') \wedge [head(e) = head(e') \vee tail(e) = tail(e')]$$

It is easy to see that if $w$ is any freely reduced word that can be read in $R_H$ along a path that starts at vertex $u$ and ends at vertex $v$, then $w$ can also be read in $\Gamma_H$ along a path that starts at vertex $\pi(u)$ and ends at vertex $\pi(v)$. Since positive words are necessarily freely reduced, it follows that $\Gamma_H$ is strongly connected. □

The converse of lemma 3.4 is also true. In fact, we prove a statement that is (a priori) stronger:

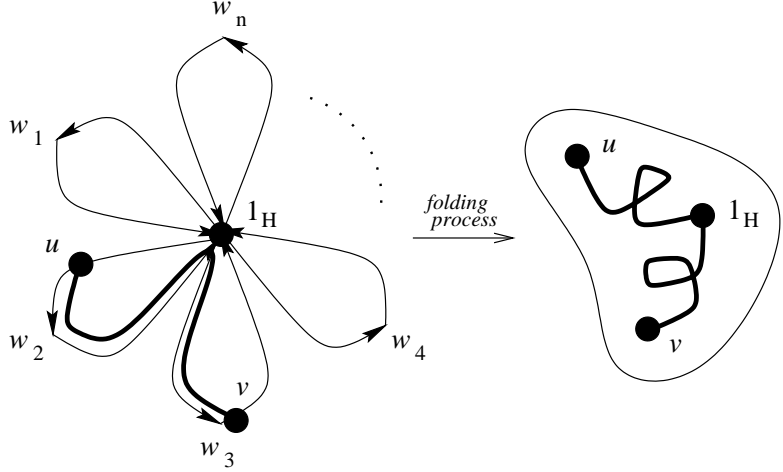

Figure 3: Positive generating set $\Rightarrow$ strong connectivity.

**Lemma 3.5.** *Let $H$ be a finitely generated subgroup of a free group $F$, and let $\Gamma = (V, E)$ be the folding of $H$. If $\Gamma_H$ is a strongly connected directed graph, then $H$ has a basis consisting of positive words.*

*Proof.* Since $\Gamma_H$ is a strongly connected directed graph, we know by corollary 3.3, that $\Gamma$ has a directed trail decomposition $P_0, P_1, \ldots, P_n$ with the property that for $i = 0, \ldots, n$, the directed graph $\Gamma_i = (V_i, E_i) = \bigcup_{j=0}^{i} P_j$ strongly connected.

Let $h_0$ be the word read along $P_0$. Now for $i = 1, \ldots, n$ put $x_i = s(P_i)$ and $y_i = t(P_i)$. Since $P_i$ is part of a *strong* directed trail decomposition, $x_i, y_i \in V_{i-1}$. Since $\Gamma_{i-1}$ is strongly connected, fix $s_i$ to be a directed trail in $\Gamma_{i-1}$ from $1_H$ to $x_i$, and fix $t_i$ to be a directed trail in $\Gamma_{i-1}$ from $y_i$ to $1_H$. Define $w_{P_i}, w_{s_i}, w_{t_i}$ to be the words read along $P_i$, $s_i$, and $t_i$ respectively. For $i = 1, \ldots, n$, put $h_i = w_{s_i} \circ w_{P_i} \circ w_{t_i}$. Clearly each $h_i$ is a positive word. Define $B_H = \{h_0, \ldots h_n\}$.

We show that $B_H$ is a free basis for $H$. For each $i = 0, \ldots, n$, fix $c_i$ be an arbitrary edge in $P_i$, and define $L_i \subset E_i$ to be the set of edges in $P_i$ excluding $c_i$. Put $T_i = (V_i, \bigcup_{j=1}^{i} L_i)$; $T_i$ is a well-defined subgraph of $\Gamma_i$. Clearly $T_0$ is spanning tree of $\Gamma_0$. Assume $T_{i-1}$ is a spanning tree of $\Gamma_{i-1}$. Then, since $P_i$ is a directed trail which attaches to $\Gamma_{i-1}$ at precisely $\{s(P_i), t(P_i)\}$, the omission of edge $c_i$ from $L_i$ suffices to ensure that $T_i$ is a spanning tree of $\Gamma_i$. By induction, $T_n$ is a spanning tree of $\Gamma_n = \Gamma_H$.

Since $B_H$ consists precisely of the Schreier transversals of $H$ relative to the spanning tree $T_n$, it follows that $B_H$ is a free basis for $H$. □

Combining the results of lemmas 3.4 and 3.5, we see that the existence of a positive generating set for a subgroup $H \leqslant_{\text{f.g.}} F$ is equivalent to the existence of a positive basis for $H$.

**Corollary 3.6.** *Let $H$ be a finitely generated subgroup of a free group $F$. $H$ is generated by a set of positive words iff $H$ has a basis consisting of positive words.*

*Proof.* $\Leftarrow$ Trivial.
$\Rightarrow$ If $H$ is generated by a set of positive words, then by lemma 3.4, $\Gamma_H$ is strongly connected. Then, by lemma 3.5, $H$ has a basis consisting of positive words. $\square$

We know from lemma 3.4 that a positively generated subgroup of (an arbitrary) free group must have a strongly connected folding. From lemma 3.2 we see that strongly connected foldings have strong directed trail decompositions. Now we show that whenever a folding of a subgroup $H \leqslant_{\text{f.g.}} F_2$ has a (not necessarily strong) directed trail decomposition, then this folding is necessarily 3-balanced.

**Lemma 3.7.** *If $H \leqslant_{\text{f.g.}} F_2$ such that $\Gamma_H = (V, E)$ has a directed trail decomposition, then $\Gamma_H$ is 3-balanced.*

*Proof.* Suppose $\Gamma_H$ has a directed trail decomposition $P_0, \ldots, P_n$. For each $i = 0, \ldots, n$, define $\Gamma_i \stackrel{\text{def}}{=} (V_i, E_i)$, where $E_i \stackrel{\text{def}}{=} \bigcup_{j=0}^{i} P_j$ and $V_i$ are the vertices induced by $E_i$. Clearly, $\Gamma_{i-1}$ and $P_i$ are subgraphs of $\Gamma_i$, and $\Gamma_n = \Gamma_H$. Note that since $(P_j)_{j=0,\ldots,i}$ is a directed trail decomposition of $\Gamma_i$, all vertices in $V_i$ have degree $\geqslant 2$.

We prove the lemma by induction on $n$.

*Base case:* When $n = 0$, $\Gamma_H$ consists of exactly one directed trail starting and ending at $1_H$. Note that then $\Gamma_H = \Gamma_0 = P_0$ cannot have vertices of odd degree. It follows that $\forall i \in \{1, 2, 3, 4\}$, $C_i(\Gamma_H) = 0$, and so the lemma holds trivially.

*Inductive step:* We will assume that the lemma holds for the folding $\Gamma_{n-1}$, and show this implies the lemma is also true for $\Gamma_n$.

Put $u = s(P_n)$ and $v = t(P_n)$. If $u, v \notin V_{n-1}$, then $u = v$; so $P_n$ (as a subgraph of $\Gamma_n$) consists only of vertices of even degree. Thus, $d_3(\Gamma_n) = d_3(\Gamma_{n-1})$, and in particular, for $i = 1, \ldots, 4$, $C_i(\Gamma_n) = C_i(\Gamma_{n-1})$. By the inductive hypothesis, the lemma holds.

By definition (2.6) of directed trail decomposition, it only remains to consider the case when *both* $u, v \in V_{n-1}$. Since $P_n$ is a trail, all vertices in $V_n \backslash V_{n-1}$ must be of even degree (either 2 or 4) in $\Gamma_n$. Since $\Gamma_{n-1}$ is a subgraph of $\Gamma_n$, it follows that $\{w \in V_{n-1} | d(w^+) > d(w)\} = \{u, v\}$. For clarity, we denote vertex $u \in V_{n-1}$ as $u^+$ when we are considering its properties as a vertex in $V_n$. For example, we denote the degree of $u$ in $\Gamma_{n-1}$ as $d(u)$, while denoting the degree of the same vertex in $\Gamma_n$ as $d(u^+)$. The terms $v, v^+, d(v)$ and $d(v^+)$ are defined analogously. Let $e_u$ [resp. $e_v$] be the first [resp. last] edge in $P_n$, and denote its label by $l_u$ [resp. $l_v$], where Note $l_u, l_v \in \{a, b\}$. The various edges, vertices and labels are depicted in figure 4.

**The case when $u = v$:** If $u = v$, then it must be that $d(u) = 2$ and $d(u^+) = 4$. So, $d_3(\Gamma_n) = d_3(\Gamma_{n-1})$, and in particular, for $i = 1, \ldots, 4$, $C_i(\Gamma_n) = C_i(\Gamma_{n-1})$. By the inductive hypothesis, the lemma holds.

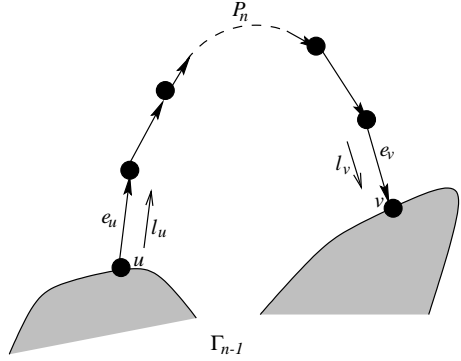

Figure 4: How $P_n$ attaches to vertices $u, v$ of $\Gamma_{n-1}$ via the edges $e_u, e_v$.

**The case when $u \neq v$:** Clearly $d(v^+) = d(v) + 1$ and $d(u^+) = d(u) + 1$. It follows that $d(u), d(v) \in \{2, 3\}$, and $d(u^+), d(v^+) \in \{3, 4\}$. We proceed now by considering each of the possible cases:

**When $d(u) = 2$:** Since $u$ has degree 2, and $\Gamma_{n-1}$ is the union of directed trails, $u$ must have one incoming and one outgoing edge in $\Gamma_{n-1}$. But since $l_u \in \{a, b\}$, it must be that (in $\Gamma_n$) $u^+$ has 2 outgoing edges and 1 incoming edge, i.e. $u^+$ is either of type $C_1$ or $C_3$. Thus, the addition of the trail $P_n$ caused the degree 2 vertex $u \in V_{n-1}$ to be transformed into a degree 3 vertex $u^+ \in V_n$ of type $C_1$ or type $C_3$. In summary, the transition $u \rightsquigarrow u^+$ causes the quantity $C_1 + C_3$ to increase by 1.

**When $d(u) = 3$:** Since vertex $u^+$ has degree 4 in $\Gamma_n$, and $l_u \in \{a, b\}$, the vertex $u^+$ has one more outgoing edge than $u$. It follows that in $\Gamma_{n-1}$, the vertex $u$ had two incoming edges and one outgoing edge, i.e. $u$ was either of type $C_2$ or $C_4$. Thus, the addition of trail $P_n$ caused the degree 3 vertex $u \in V_{n-1}$ whose type was either $C_2$ or $C_4$ to be transformed into a degree 4 vertex $u^+ \in V_n$. In summary, the transition $u \rightsquigarrow u^+$ causes the quantity $C_2 + C_4$ to decrease by 1.

**When $d(v) = 2$:** Since $v$ has degree 2, and $\Gamma_{n-1}$ is the union of directed trails, $v$ must have one incoming and one outgoing edge in $\Gamma_{n-1}$. But since $l_v \in \{a, b\}$, it must be that (in $\Gamma_n$) $v^+$ has 2 incoming edges and 1 outgoing edge, i.e. $v^+$ is either of type $C_2$ or $C_4$. Thus, the addition of the trail $P_n$ caused the degree 2 vertex $v \in V_{n-1}$ to be transformed into a degree 3 vertex $v^+ \in V_n$ of type $C_2$ or type $C_4$. In summary, the transition $v \rightsquigarrow v^+$ causes the quantity $C_2 + C_4$ to increase by 1.

**When $d(v) = 3$:** Since vertex $v^+$ has degree 4 in $\Gamma_n$, and $l_v \in \{a, b\}$, the vertex $v^+$ has one more incoming edge than $v$. It follows that in $\Gamma_{n-1}$, the vertex $v$ had two outgoing edges and one incoming edge, i.e. $v$ was either of type $C_1$ or type $C_3$. Thus, the addition of trail $P_n$ caused the degree 3 vertex $v \in V_{n-1}$ whose type was either $C_1$ or $C_3$ to be transformed into a degree 4 vertex $v^+ \in V_n$. In summary, the transition $v \rightsquigarrow v^+$ causes the quantity $C_1 + C_3$ to decrease by 1.

The conclusions of these arguments are tabulated below.

|  | $d(v) = 2$ | $d(v) = 3$ |
|---|---|---|
| $d(u) = 2$ | $v \rightsquigarrow v^+$: $C_2 + C_4$ increases by 1, $u \rightsquigarrow u^+$: $C_1 + C_3$ increases by 1 | $v \rightsquigarrow v^+$: $C_1 + C_3$ decreases by 1, $u \rightsquigarrow u^+$: $C_1 + C_3$ increases by 1 |
| $d(u) = 3$ | $v \rightsquigarrow v^+$: $C_2 + C_4$ increases by 1, $u \rightsquigarrow u^+$: $C_2 + C_4$ decreases by 1 | $v \rightsquigarrow v^+$: $C_1 + C_3$ decreases by 1, $u \rightsquigarrow u^+$: $C_2 + C_4$ decreases by 1 |

Assuming the induction hypothesis that $C_1(\Gamma_{n-1}) + C_3(\Gamma_{n-1}) = C_2(\Gamma_{n-1}) + C_4(\Gamma_{n-1})$ the table above shows that in each of the 4 possible cases, $C_1(\Gamma_n) + C_3(\Gamma_n) = C_2(\Gamma_n) + C_4(\Gamma_n)$. By induction then, the lemma holds. □

The following remark will be used to argue that a subsequent theorem about subgroups of $F_2$ holds in general for the finitely generated subgroups of any free group $F$.

**Remark 3.8.** Let $\phi_n$ be the homomorphism of $F_n = F(\{x_1, \ldots, x_n\})$ into $F_2 = F(\{a, b\})$ defined by $\phi_n : x_i \mapsto a^i b a^i$, $(i = 1, \ldots, n)$. It is easy to verify that $\phi_n$ is an embedding which takes positive words in $F_n$ to positive words in $F_2$. In particular, if $H \leqslant_{\text{f.g.}} F_n$, then $rank\ \phi(H) = rank\ H$, and if $H$ has a positive generating set, then $\phi(H)$ has a positive generating set.

The main theorem may now be proved:

**Theorem 3.9.** *If $H, K$ are two finitely generated subgroups of a free group $F$ and at least one of the two subgroups is generated by a set of positive words, then the pair $(H, K)$ satisfy the Hanna Neumann conjecture.*

*Proof.* WLOG, let $H$ have a positive generating set $B_H = \{h_1, \cdots h_n\}$.

Suppose first that $F = F_2$. By lemma 3.4, this implies $\Gamma_H$ is strongly connected. By lemma 3.2, a directed trail decomposition of $\Gamma_H$ exists. By lemma 3.7, $\Gamma_H$ is 3-balanced. Finally, by remark 2.2, $H$ cannot be part of any counterexample to the conjecture. This proves the case when $F = F_2$.

Now suppose $F \neq F_2$. Since $H, K$ are finitely generated, WLOG, we can assume that $F = F_n$ for some finite $n$. If $H, K$ were a counterexample to the conjecture, then by remark 3.8, $\phi_n(H), \phi_n(K) \leqslant_{\text{f.g.}} F_2$ and $\phi_n(H)$ has a positive generating set. Moreover, $rank\ H = rank\ \phi_n(H)$, $rank\ K = rank\ \phi_n(K)$, $rank\ H \cap K = rank\ \phi_n(H \cap K) = rank\ (\phi_n(H) \cap \phi_n(K))$. Thus $\phi_n(H), \phi_n(K)$ are a counterexample to the conjecture, contradicting our proof above for the case when $F = F_2$. This proves the case when $F \neq F_2$. □

**Corollary 3.10.** *If $H, K$ are two finitely generated subgroups of a free group $F$ and at least one of the two subgroups is potentially positive, then the pair $(H, K)$ satisfy the Hanna Neumann conjecture.*

*Proof.* Clearly, if $H, K \leqslant_{\text{f.g.}} F$ are a counterexample to the conjecture and $\phi \in Aut(F)$, then $\phi(H), \phi(K)$ are also a counterexample to the conjecture. Thus theorem 3.9 continues to hold when the condition of one group being positively generated is replaced with the condition of one group being potentially positive. □

## 4 Further Analysis

The proof of theorem 3.9 hinges on the observation that finitely generated subgroups of $F_2$ which have 3-balanced foldings cannot be part of any counterexample to the conjecture, and foldings which have directed trail decompositions are necessarily 3-balanced. The success of the approach naturally leads us to inquire about sufficient conditions for a folding to possess a directed trail decomposition. The next lemma answers this question for finitely generated subgroups of $F_2$.

**Lemma 4.1.** *Let $H \leqslant_{f.g.} F_2$ and $\Gamma = (V, E)$ be the folding of $H$. Then $\Gamma$ has no sources and no sinks if and only if $\Gamma$ has a directed trail decomposition.*

*Proof.* ⇒ Decompose $\Gamma$ into strongly connected components. Let $\widetilde{J}_1, \ldots \widetilde{J}_m$ denote those strongly connected components whose size (number of vertices) is $> 1$. Take

$$G_0 = \bigcup_{i=1,\ldots,m} \widetilde{J}_i$$

Since each $\widetilde{J}_i$ is strongly connected, by lemma 3.2, each $\widetilde{J}_i$ has a directed trail decomposition $Q_i$. Since the $\widetilde{J}_i$ are pairwise disjoint, it follows that $G_0$ has a directed trail decomposition $P_0$—simply take $P_0$ to be $Q_1, \ldots, Q_m$. We extend $P_0$ to a directed trail decomposition of $\Gamma$ in stages. At each successive stage $i$ (starting at $i = 0$):

1. If $E[\Gamma] \backslash E[G_i] = \emptyset$, halt. Otherwise, select any directed edge $(u, v) \in E[\Gamma] \backslash E[G_i]$ (from $u$ to $v$). Starting at vertex $u$ we walk backwards along (arbitrarily chosen) incoming edges, until reach a vertex $u' \in V[G_i]$. Likewise, starting at vertex $v$ we walk forwards along (arbitrarily chosen) outgoing edges until we reach a vertex $v' \in V[G_i]$. We cannot get stuck in either of these steps, because $\Gamma$ contains neither sources nor sinks; we cannot get trapped in a loop before we find a vertex in $V[G_i]$ because then we have discovered a strongly connected component that must have been omitted from the set $\widetilde{J}_1, \ldots \widetilde{J}_m$, a contradiction.

   Define $P_{i+1}$ to be the directed trail from $u' \rightsquigarrow u \rightsquigarrow v \rightsquigarrow v'$ described above. Notice that $P_{i+1}$ attaches $G_i$ at precisely its endpoints $u', v'$. We append $P_{i+1}$ to the directed trail decomposition at stage $i$, obtaining a trail decomposition of $G_{i+1} = G_i \cup P_{i+1}$. Increment $i$, then repeat step 1.

At the end of this procedure, we have constructed a directed trail decomposition of $\Gamma$, as claimed.

⇐ Suppose that $\Gamma = (V, E)$ has a directed trail decomposition $P_0, \ldots, P_k$ and (towards contradiction) also has a source vertex $v_0 \in V$. Since $\{P_i\}_{i=0,\ldots,k}$ covers $E[\Gamma]$, let $i_0 \in \{0, 1, \ldots, k\}$ be the least integer for which $v_0 \in V[P_{i_0}]$.

Suppose $v_0 \neq s(P_{i_0}), t(P_{i_0})$; then $v_0$ has two outgoing edges, contradicting the fact that $P_{i_0}$ is a directed trail.

Suppose $v_0 \in \{s(P_{i_0}), t(P_{i_0})\}$; then minimality of $i_0$ implies that $V[P_{i_0}] \cap \bigcup_{j=0}^{i_0-1} V[P_j] = \emptyset$. So, by definition 2.6 of directed trail decomposition, the endpoints of $P_{i_0}$ must coincide. But $P_{i_0}$ is not a directed trail, since its final edge is not oriented towards $t(P_{i_0})$.

A completely analogous argument shows if $\Gamma$ has a sink vertex, then $\Gamma$ cannot have a directed trail decomposition. □

It is well-known that if the Hanna Neumann conjecture holds for subgroups of $F_2$, then it holds in general for subgroups of any free group $F$. Lemma 4.1 then has the following consequence for the conjecture:

**Theorem 4.2.** *If $H, K$ are two finitely generated subgroups of the free group $F_2$ and at least one of the two subgroups is source/sink-free, then the pair $(H, K)$ satisfy the Hanna Neumann conjecture.*

*Proof.* If $\Gamma_H$ has neither source nor sink vertices, then by lemma 4.1, $\Gamma_H$ has a directed trail decomposition. So, lemma 3.7 applies and hence $\Gamma_H$ must be 3-balanced. Then, by remark 2.2, $H$ cannot be part of any counterexample to the conjecture. □

**Corollary 4.3.** *If $H, K$ are two finitely generated subgroups of the free group $F_2$ and at least one of the two subgroups is potentially source/sink-free, then the pair $(H, K)$ satisfy the Hanna Neumann conjecture.*

*Proof.* Clearly, if $H, K \leqslant_{\text{f.g.}} F_2$ are a counterexample to the conjecture and $\phi \in Aut(F_2)$, then $\phi(H), \phi(K)$ are also a counterexample to the conjecture. Thus theorem 4.2 continues to hold when the condition of one group being source/sink-free is replaced with the condition of one group being potentially source/sink-free. □

# 5 Acknowledgements

I would like to thank the weekly Group Theory Seminar at the City University of New York Graduate Center, where I first came to know of these problems. My gratitude to Alexei Miasnikov and Toshiaki Jutsikawa for many invaluable suggestions and spirited discussions along the way.

I gratefully acknowledge the Mathematics department at the City University of New York Graduate Center, for funding this work as part of the my ongoing doctoral research. I also thank both Advanced Engineering & Sciences, and the Center for Computational Sciences at the Naval Research Laboratory in Washington DC for their support of these endeavors.

Bilal Khan
Department of Mathematics, City University of New York Graduate Center,
365 5th Avenue, NY, NY 10016.